# Polyhedra of Genus 2 with 10 Vertices and Minimal Coordinates.

28th July 2005

### Authors

Stefan Hougardy, Frank H. Lutz, and Mariano Zelke

### Description

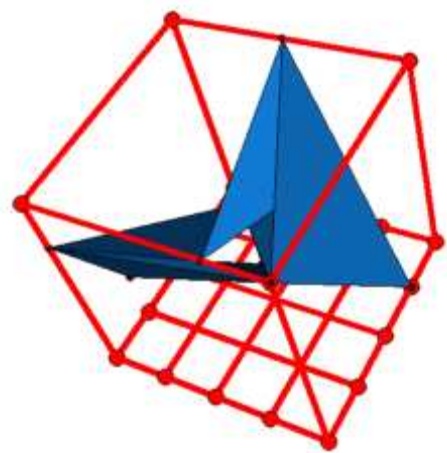

We give coordinate-minimal geometric realizations in general position of all 865 vertex-minimal triangulations of the orientable surface of genus 2 in the 4x4x4-cube.

In 1890, Heawood [8] proved that a triangulation of a (closed) surface M of Euler characteristic chi(M) has at least

$$n \geq 1/2(7+\sqrt{49-24 \cdot \text{chi}(M)})$$

vertices. With the exception of the orientable surface of genus 2, the Klein bottle, and the non-orientable surface of genus 3 this bound is tight, as was shown only much later by Ringel [15] for non-orientable surfaces and by Jungerman and Ringel [11] for orientable surfaces. In the three exceptional cases one extra vertex has to be added, respectively, to allow for triangulations. For the orientable surface of genus 2 this fact was proved by Huneke [10], thus establishing that vertex-minimal triangulations of the orientable surface of genus 2 have 10 vertices.

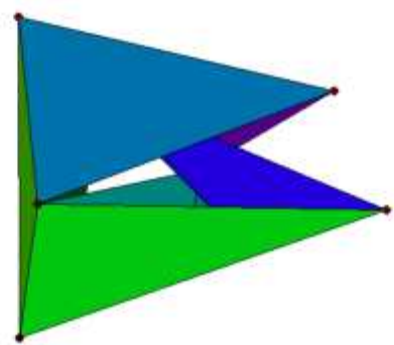

A complete enumeration of all 42426 triangulated surfaces with 10 vertices was obtained in [12]; see [13] for a list of facets of the triangulations. In particular, there are exactly 865 combinatorially distinct vertex-minimal 10-vertex triangulations of the orientable surface of genus 2.

By Steinitz' theorem (cf. [16, Ch. 4]), every triangulated 2-sphere is realizable geometrically as the boundary complex of a convex 3-dimensional polytope. For triangulations of orientable surfaces of genus $g \geq 1$ it was asked by Grünbaum [7, Ch. 13.2] whether they can always berealized geometrically as a polyhedron in $R^3$, i.e., with straight edges, flat triangles, and without self intersections? In general, the answer turned out to be NO: Bokowski and Guedes de Oliveira [3] showed that there are non-realizable triangulations for surfaces of larger genus, $g \geq 6$. However, for surfaces of genus $1 \leq g \leq 5$ the problem remains open.

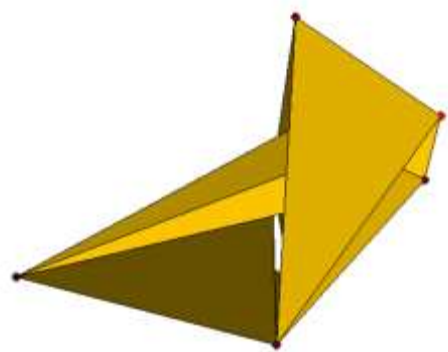

Geometric realizations for several examples of triangulated orientable surfaces of genus 2, 3, and 4 with respective minimal numbers of vertices 10, 10, and 11 were constructed by Brehm and Bokowski [1], [2], [5], [6].

In fact, the realizability problem for triangulated surfaces is decidable (cf. Bokowski and Sturmfels [4, p. 50]), but there is no algorithm known that would solve the realization problem for instances with, say, 10 vertices in reasonable time.

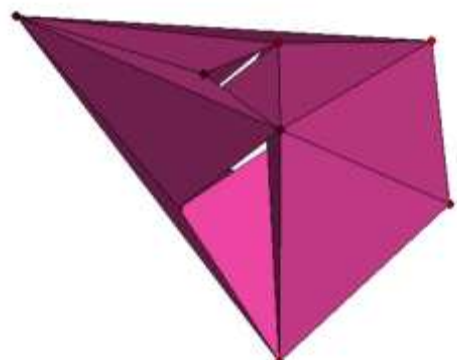

By random realization, geometric realizations in the 32768x32768x32768-cube were obtained by Lutz [12] for 864 of the 865 examples of vertex-minimal 10-vertex triangulations of the orientable surface of genus 2. A realization of the remaining example was constructed by Bokowski.

Theorem (Bokowski and Lutz, cf. [12]): All 865 vertex-minimal 10-vertex triangulations of the orientable surface of genus 2 can be realized geometrically in $R^3$.

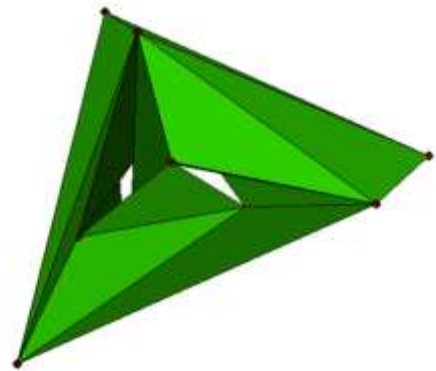

Surprisingly, all these examples have realizations with small coordinates.

Theorem: All 865 vertex-minimal 10-vertex triangulations of the orientable surface of genus 2 have realizations in general position in the 4x4x4-cube, but cannot be realized in general position in the 3x3x3-cube.

Our realization algorithm for general position realizations with small coordinates is a variant of the isomorph-free exhaustive generation as described by McKay [14] for classes of objects with an inductive construction process. In our case, we generate sets of increasing size s with up to 10 points with integer coordinates in general position in the nxnxn-cube, which allow realizations of at least one induced subcomplex with s vertices of a given triangulated orientable surface.

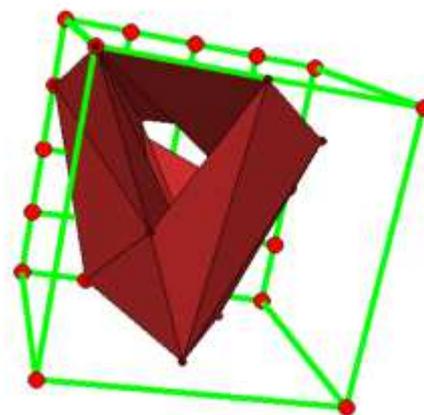

This realization algorithm can easily be adapted to also obtain proper realizations, i.e., realizations that do not have coplanar neighboring triangles, but which not necessarily need to be in general position (and which therefore might be of even smaller size). Still more general, realizations with coplanar neighboring triangles can be produced as well. This way, for example, we found realizations of triangulations of the torus in the 2x2x2-cube; see [9].

Remark: The six displayed examples were selected from the collection of the 865 general position realizations in the 4x4x4-cube for their clearly visible holes.

**Authors' Addresses**

Stefan Hougardy

> Humboldt-Universität zu Berlin
> Institut für Informatik
> Unter den Linden 6
> 10099 Berlin
> Germany
> hougardy@informatik.hu-berlin.de
> http://www.informatik.hu-berlin.de/~hougardy/

Frank H. Lutz

> Technische Universität Berlin
> Fakultät II - Mathematik und Naturwissenschaften
> Institut für Mathematik, Sekr. MA 6-2
> Straße des 17. Juni 136
> 10623 Berlin
> Germany
> lutz@math.tu-berlin.de
> http://www.math.tu-berlin.de/~lutz

Mariano Zelke

> Humboldt-Universität zu Berlin
> Institut für Informatik
> Unter den Linden 6
> 10099 Berlin
> Germany
> zelke@informatik.hu-berlin.de